\documentclass[11pt,final]{amsart}

\usepackage[
  a4paper,
  margin=2.5cm
]{geometry}
\usepackage[T1]{fontenc}
\usepackage[utf8]{inputenc}
\usepackage{lmodern}
\usepackage{microtype}
\usepackage{mathtools}
\usepackage{amssymb}
\usepackage{mathrsfs}
\usepackage{bm}
\usepackage{enumitem}
\usepackage{booktabs}
\usepackage{xcolor}
\usepackage{hyperref}
\usepackage{aliascnt}
\usepackage{tikz}
\usepackage[nameinlink,capitalize,noabbrev]{cleveref}
\usetikzlibrary{arrows.meta, intersections}

\hypersetup{
colorlinks=true,
linkcolor=black
}

\numberwithin{equation}{section}

\newcommand{\dd}{\,\mathrm{d}}

\newtheorem{theorem}{Theorem}[section]
\newaliascnt{proposition}{theorem}
\newtheorem{proposition}[proposition]{Proposition}
\aliascntresetthe{proposition}
\newaliascnt{lemma}{theorem}
\newtheorem{lemma}[lemma]{Lemma}
\aliascntresetthe{lemma}
\newaliascnt{corollary}{theorem}

\aliascntresetthe{corollary}
\theoremstyle{definition}
\newaliascnt{definition}{theorem}

\aliascntresetthe{definition}
\theoremstyle{remark}
\newaliascnt{remark}{theorem}

\aliascntresetthe{remark}
\title[Zeros Counting of MIT Bag Conditions]{Upper Bounds for the Number of Zeros of Dirac Eigenspinors on Planar Domains with Infinite Mass Boundary Conditions}

\author{Zhang Yitian}

\date{}

\begin{document}

\begin{abstract}
We study the zeros, counted with multiplicity, of eigenspinors of the massless Dirac operator with MIT bag boundary conditions on smooth, bounded, connected planar domains. We get a strict Pleijel-type asymptotic bound, a zero-counting inequality in terms of the Morse index at sufficiently high frequencies, and an upper bound of the form $k+C_\Omega$.
\end{abstract}

\maketitle

\section{Introduction}\label{sec:introduction}

The distribution of zeros of eigenfunctions is a classical subject in spectral theory. For one--dimensional Sturm--Liouville problems with separated boundary conditions, the Sturm--Liouville theorem provides an identity of the number of zeros of eigenfunctions: an eigenfunction associated with the $k$-th eigenvalue has exactly $k-1$ zeros in the interior of the interval. In higher dimensions, the geometry of the zero sets becomes more intricate, and the focus shifts to the nodal sets of eigenfunctions. 

Let $\Omega\subset\mathbb{R}^n$ be a bounded domain, and consider the Dirichlet eigenvalue problem for the Laplacian:
\[
\begin{cases}
-\Delta u_k = \mu_k u_k,
& x\in\Omega,\\
u_k=0,
& x\in\partial\Omega,
\end{cases}
\]
where, counted with multiplicity:
\[
0<\mu_1\leqslant\mu_2\leqslant\cdots.
\]
Denote the nodal set of $u_k$ by
\[
Z(u_k)\coloneqq \{x\in\Omega:u_k(x)=0\}.
\]
Courant's nodal domain theorem \cite{Courant1923} asserts that
\[
\#\pi_0(\Omega\setminus Z(u_k))\leqslant k,
\]
that is, the $k$-th eigenfunction has at most $k$ nodal domains. Pleijel's nodal domain theorem \cite{Pleijel} asserts that, when $\Omega\subset\mathbb{R}^2$,
\[
\limsup_{k\to\infty}\frac{\#\pi_0(\Omega\setminus Z(u_k))}{k}
\leqslant
\frac{4}{j_{0,1}^2}
\approx 0.69166<1,
\]
where $j_{0,1}$ denotes the first positive zero of the Bessel function $J_0$.

But Courant-type estimates do not hold for all spectral problems. Recently, Enciso--Pistoia--Provenzano \cite{EPP2026} proved that, there is no universal upper bound on the number of boundary nodal domains of Steklov eigenfunctions depending only on the spectral index. This result illustrates that the geometry of zero sets may differ fundamentally in different elliptic spectral problems.

Unlike the Laplacian, the Dirac operator is a first-order elliptic system whose eigenfunctions are complex vector-valued spinors. Consequently, there is generally no natural decomposition into positive and negative nodal domains. On a planar domain $\Omega\subset\mathbb{R}^2$ (or $\mathbb{C}$), the two-dimensional massless Dirac operator is defined by
\[
D\coloneqq -i\sigma\cdot\nabla
=-i(\sigma_1\partial_x+\sigma_2\partial_y),
\]
where $\sigma=(\sigma_1,\sigma_2)$ and
$\sigma_1,\sigma_2$ are the Pauli matrices
\[
\sigma_1=
\begin{pmatrix}
0&1\\
1&0
\end{pmatrix},
\qquad
\sigma_2=
\begin{pmatrix}
0&-i\\
i&0
\end{pmatrix}.
\]
By the Wirtinger derivatives
\[
\partial_z= \frac{1}{2}(\partial_x-i\partial_y),
\qquad
\partial_{\bar{z}}
= \frac{1}{2}(\partial_x+i\partial_y),
\]
we may express the Dirac operator as
\[
D=-2i
\begin{pmatrix}
0&\partial_z\\
\partial_{\bar{z}}&0
\end{pmatrix}.
\]
For a spinor $\psi\coloneqq (u,v)^{\mathsf{T}}$, its zero set is defined as the common zero set of its two complex components: $Z(\psi)\coloneqq \{x\in\Omega:u(x)=v(x)=0\}$.

B\"ar \cite{Baer} established a structural theorem for the zero sets of solutions to Dirac equations. In particular, in two-dimensional case, the zeros of any nontrivial solution are isolated and have finite vanishing order. Thus, in contrast to the number of nodal domains for scalar eigenfunctions, the sum of the vanishing orders over all common zeros provides a natural counting quantity for studying the distribution of zeros of two-dimensional Dirac eigenspinors. For a closed two-dimensional Riemannian spin manifold $(M^2,g)$, Branding \cite{Branding} established, using the Bochner method, an inequality between the total vanishing order of a Dirac eigenspinor and its corresponding eigenvalue. More precisely, if $D_M\psi=\lambda\psi$, then
\begin{equation}\label{eq:Branding-Bochner}
N(\psi)
\leqslant
\frac{\operatorname{Area}(M)}{4\pi}\lambda^2-\frac{\chi(M)}{2},
\end{equation}
where
\[
N(\psi)\coloneqq \sum_{p\in Z(\psi)}\operatorname{ord}_p\psi
\]
denotes the total vanishing order of the spinor at its zeros. $\operatorname{ord}_p\psi$ denotes the vanishing order of $\psi$ at a zero $p$ and $\chi(M)$ is the Euler characteristic of $M$. Moreover, The Dirac operator on a closed spin surface satisfies Weyl's asymptotic law. Indeed, the Lichnerowicz formula
\[
D_M^2
=
\nabla^{\Sigma *}\nabla^\Sigma+\frac K2
\]
shows that $D_M^2$ is an operator of Laplace type, where $K$ denotes the Gaussian curvature. By the general Weyl law for elliptic operators \cite{Hormander1968,Ivrii}, we have
\[
N_{D_M^2}(T)
=
\frac{\operatorname{Area}(M)}{2\pi}T+o(T),
\qquad T\to\infty,
\]
where $N_{D_M^2}(T)$ denotes the number of eigenvalues of $D_M^2$ strictly less than $T$, counted with multiplicity. Owing to the chiral symmetry of the Dirac operator on a closed two-dimensional spin manifold, its nonzero spectrum is symmetric about the origin, with $\lambda$ and $-\lambda$ having the same multiplicity \cite{LawsonMichelsohn1989}. Thus, listing the positive eigenvalues in nondecreasing order, counted with complex multiplicity,
\[
0<\lambda_1^+\leqslant\lambda_2^+\leqslant\cdots,
\]
we obtain
\begin{equation}\label{eq:M_2-lambda-weyl}
\frac{\operatorname{Area}(M)}{4\pi}
(\lambda_k^+)^2
= k+o(k),
\qquad 
k\to\infty.
\end{equation}
Combining \eqref{eq:Branding-Bochner} with \eqref{eq:M_2-lambda-weyl}, we obtain
\[
N(\psi_k)\leqslant k+o(k)
\]
for any eigenspinor $\psi_k$ associated with $\lambda_k^+$. This shows that, on closed two-dimensional spin manifolds, the total vanishing order of Dirac eigenspinors admits an asymptotically linear upper bound in terms of the spectral index. 

In this paper, we consider the corresponding problem on planar domains with boundary. Let $\Omega\subset\mathbb{C}$ be a bounded connected domain with $C^\infty$ boundary $\partial\Omega$. We denote by $\nu=(\nu_1,\nu_2)$ the outward unit normal and by $\tau=(\tau_1,\tau_2)\coloneqq(-\nu_2,\nu_1)$ the positively oriented unit tangent. Define
\[
\sigma_\tau
\coloneqq
\tau_1\sigma_1+\tau_2\sigma_2
=
\begin{pmatrix}
0&\tau_1-i\tau_2\\
\tau_1+i\tau_2&0
\end{pmatrix}.
\]
We consider the eigenvalue problem for the two-dimensional massless Dirac operator with infinite mass boundary conditions:
\begin{equation}\label{eq:mit-model}
\begin{cases}
D\psi= \lambda\psi,
& x\in\Omega,\\
\sigma_\tau\psi= \psi,
& x\in\partial\Omega.
\end{cases}
\end{equation}
This boundary condition is also known as the MIT bag boundary condition (the corresponding operator is called MIT bag Dirac operator). $\psi$ is called the eigenspinor of MIT bag Dirac operator. For $\lambda\neq 0$, the eigenvalue problem \eqref{eq:mit-model} is equivalent to the following $\bar\partial$-Robin Laplacian problem:
\begin{equation}\label{eq:cn-robin-problem}
\begin{cases}
-\Delta u= \lambda^2u,
& x\in\Omega,\\
2\bar\nu\partial_{\bar{z}}u+\lambda u= 0,
& x\in\partial\Omega.
\end{cases}
\end{equation}
Benguria--Fournais--Stockmeyer--Van Den Bosch \cite{BFSV} established self-adjointness for a class of boundary conditions, including the infinite mass boundary condition. Moreover, the MIT bag Dirac operator on a smooth bounded domain has compact resolvent, and its spectrum consists entirely of real eigenvalues of finite multiplicity.

Duran and Duran--Mas--Sanz-Perela \cite{Duran,DMS} studied the variational properties of the $\bar\partial$-Robin Laplacian and its spectral correspondence with quantum-dot Dirac operators, thereby providing a new approach to Dirac spectral problems via quadratic forms.

More recently, Ammari--Qiu \cite[Section 4]{AQ} investigated the structure of boundary zeros of eigenspinors by means of local conformal flattening, spinor reflection, and holomorphic factorization. 

Since boundary zeros also have finite integer vanishing orders and their contributions to the local flux integrals are half those of interior zeros of the same order, we define the weighted total vanishing order by
\begin{equation}\label{eq:cn-zero-count}
\mathfrak{N}(\psi)
\coloneqq
\sum_{p\in Z(\psi)\cap\Omega}
\operatorname{ord}_p\psi
+
\frac{1}{2}
\sum_{p\in Z(\psi)\cap\partial\Omega}
\operatorname{ord}_p\psi.
\end{equation}
where here
\[
Z(\psi)\coloneqq \left\{x\in\overline\Omega:\psi(x)=0\right\}.
\]

Our aim is to establish quantitative estimates relating the weighted total vanishing order of eigenspinors of \eqref{eq:mit-model} to the index of the corresponding positive eigenvalue. Our approach is based on extending Branding's Bochner method to the setting of manifolds with boundary. More precisely, we exploit the nonnegative defect term in an exact Bochner identity, together with a Morse index comparison for the $\bar\partial$-Robin quadratic form, to control the weighted total vanishing order of Dirac eigenspinors. More precisely, let
\[
0<\lambda_1^+\leqslant\lambda_2^+\leqslant\cdots
\]
be the positive eigenvalues of the MIT Dirac operator,
listed in nondecreasing order and counted with multiplicity,
and let $\psi_k$ be a corresponding nonzero eigenspinor.
We establish the estimate
\[
\mathfrak{N}(\psi_k)
\leqslant k+C_\Omega.
\]

Section~\ref{sec:preliminaries} introduces the necessary preliminaries; Section~\ref{sec:defect_identity} gives the exact Bochner identity; In Section~\ref{sec:lower_bound_defect}, we prove the positive $L^1$ lower bound of the defect; In Section~\ref{sec:Pleijel-type-thm}, we prove the Pleijel--type theorem of eigenspinors; In Section~\ref{sec:morse-index}, we derive the $k+C_\Omega$ bound.

\section{Preliminaries}
\label{sec:preliminaries}

In this section, we recall the basic spectral properties of the two-dimensional MIT bag Dirac operator, the structure of the zero sets of its eigenspinors, and the variational connection between the MIT bag Dirac operator and the $\bar\partial$-Robin Laplacian.

Let $\psi=(u,v)^{\mathsf T}\not\equiv 0$ be an eigenspinor of $D$. By B\"ar's structural theorem for zero sets \cite{Baer}, together with the local structure of the
Dirac equation and Taylor expansion, we obtain the following local representation near an interior zero $p$. Choosing a complex coordinate $z=x+iy$ centered at $p$, there exist an integer $m\geqslant 1$ and a nonzero vector $c\in\mathbb{C}^2$ such that
\begin{equation}\label{eq:ord-zero}
(u,\bar{v})= z^m c+R(z),
\qquad
R(z)= (r^{m+1}),\quad
\nabla R(z)= O(r^m),
\end{equation}
where $r=|z|$. For common zeros on the boundary, Ammari--Qiu \cite{AQ} established a finite-order local expansion. More precisely, in suitable conformal coordinates on
a half-disk, the two components of the transformed spinor have the same finite integer vanishing order. Consequently, all zeros of a nontrivial eigenspinor in $\overline\Omega$ are isolated, and hence the zero set is finite. Thus, whether $p$ is an interior zero or a common zero on the boundary, we have, in the corresponding local coordinates,
\[
|\psi|
=
r^m\left(c_p+O(r)\right),
\qquad c_p>0,
\]
and
\[
\partial_r\log|\psi|
=
\frac{m}{r}+O(1),
\qquad
m=\operatorname{ord}_p\psi.
\]
Consequently, upon removing a small disk or half-disk of radius $\varepsilon$ around each zero and taking $\nu_\varepsilon$ to be the outward unit normal to the remaining domain, the boundary integral around an interior zero satisfies
\begin{equation}\label{eq:interior-zero-boundary-integral}
\int_{\partial B_\varepsilon(p)}
\partial_{\nu_\varepsilon}\log|\psi|\dd s
=
-2\pi m+o(1),
\end{equation}
whereas the corresponding integral along the circular arc at a boundary zero satisfies
\begin{equation}\label{eq:boundary-zero-boundary-integral}
\int_{\partial B_\varepsilon(p)\cap\Omega}
\partial_{\nu_\varepsilon}\log|\psi|\dd s
=
-\pi m+o(1),
\end{equation}
as $\varepsilon\to 0^+$. It follows that the contribution of a boundary zero is exactly half that of an interior zero of the same vanishing order. Accordingly, we define the weighted total vanishing order of an eigenspinor by \eqref{eq:cn-zero-count}. A similar convention of counting boundary zeros with half weight also appears in \cite[Appendix]{Wendl2010}. 

We now turn to the equivalent scalar problem \eqref{eq:cn-robin-problem}. Following the variational approach of Antunes--Benguria--Lotoreichik--Ourmi\`eres-Bonafos \cite{Antunes2021}, for $t>0$ we define the quadratic form
\begin{equation}
\label{eq:cn-form}
Q_t[w]
\coloneqq
4\int_\Omega|\partial_{\bar{z}}w|^2\dd x
+t\int_{\partial\Omega}|\operatorname{tr} w|^2\dd s
-t^2\int_\Omega|w|^2\dd x.
\end{equation}
The form domain is given by
\[
E(\Omega)
\coloneqq
\left\{
w\in L^2(\Omega):
\partial_{\bar{z}}w\in L^2(\Omega),\
\operatorname{tr} w\in L^2(\partial\Omega)
\right\},
\]
where $\partial_{\bar{z}}w$ is understood in the distributional sense, and $\operatorname{tr}w$ denotes the trace of $w$. Antunes--Benguria--Lotoreichik--Ourmi\`eres-Bonafos \cite{Antunes2021} constructed, by completion, the densely defined and closed quadratic form $Q_t$, characterized its form domain $E(\Omega)$, and established the continuous embedding
\[
E(\Omega)\hookrightarrow H^{1/2}(\Omega).
\]
In particular, the Sobolev compact embedding theorem yields
\[
E(\Omega)
\hookrightarrow\hookrightarrow L^2(\Omega).
\]
For the Hilbert space structure of $E(\Omega)$ and its compact embedding into $L^2(\Omega)$, we also refer to Duran \cite[Lemma 3.1]{Duran}, while the closedness of
the corresponding $\bar\partial$-Robin quadratic form is established in \cite[Lemma 3.2]{Duran}. Consequently, for each $t>0$, there exists a unique self-adjoint operator $A_t$ associated with $Q_t$, and $A_t$ has compact resolvent \cite[Corollary 30]{Antunes2021}. We denote by $\operatorname{ind}Q_t$ the Morse index of $Q_t$, defined as the maximal complex dimension of a subspace on which $Q_t$ is negative definite, namely,
\[
\operatorname{ind}Q_t
\coloneqq
\sup\left\{
\dim_{\mathbb{C}}L:
L\subset E(\Omega),\
Q_t[w]<0,\ 
\forall w\in L\setminus\{0\}
\right\}.
\]
Since $A_t$ is self-adjoint, bounded from below, and has compact resolvent, $\operatorname{ind}Q_t$ equals the number of negative eigenvalues of $A_t$, counted with
multiplicity. 

Finally, Duran--Mas--Sanz-Perela \cite{DMS} established a one-to-one correspondence between the solutions of the quantum-dot Dirac eigenvalue problem and those of the corresponding $\bar\partial$-Robin eigenvalue problem. In particular, for the massless Dirac operator with MIT bag boundary conditions considered in this paper, the mapping
\[
u\longmapsto
\begin{pmatrix}
u\\
-\dfrac{2i}{t}\partial_{\bar{z}}u
\end{pmatrix}
\]
induces, for every $t>0$, a complex-linear isomorphism
\begin{equation}\label{eq:cn-kernel-isomorphism}
\ker A_t\cong\ker(D-tI).
\end{equation}

Now we give the main result of this paper.
\begin{theorem}\label{thm:cn-main}
Denote $E_\lambda\coloneqq \ker(D-\lambda)$. There exists $\eta_*>0$, such that
\begin{equation}\label{eq:cn-main-pleijel}
    \limsup_{k\to\infty}
    \sup_{\psi\in E_{\lambda_k^+}\setminus\{0\}}
    \frac{\mathfrak{N}(\psi)}{k}
    \leqslant 1-\eta_*<1.
\end{equation}
Moreover, there exists $L_\Omega$ depending only on $\Omega$, such that $\forall\lambda\geqslant L_\Omega$, we have
\begin{equation}\label{eq:cn-main-high}
\mathfrak{N}(\psi)\leqslant \operatorname{ind} Q_\lambda.
\end{equation}
Then there exists $C_\Omega$ depending only on $\Omega$, such that $\forall \lambda> 0$, we have
\begin{equation}\label{eq:cn-main-all}
\mathfrak{N}(\psi)\leqslant \operatorname{ind} Q_\lambda+ C_\Omega.
\end{equation}
Thus, $\forall \psi\in E_{\lambda_k^+}\setminus\{0\}$,
\begin{equation}\label{eq:cn-main-courant}
    \mathfrak{N}(\psi)\leqslant k-1+ C_\Omega.
\end{equation}
\end{theorem}

\section{A Bochner Identity}\label{sec:defect_identity}
Let $s$ denote the arc-length parameter along the direction of $\tau$, then 
\begin{equation}\label{eq:frenet-convention}
\partial_s\tau=-\kappa\nu,
\qquad
\partial_s\nu=\kappa\tau.
\end{equation}
Along the boundary, the Dirac operator can thus be decomposed as
\begin{equation}\label{eq:dirac-boundary-decomposition}
D= -i\sigma_\nu\partial_\nu-i\sigma_\tau\partial_s.
\end{equation}

By the Clifford anticommutation relations for the Pauli matrices, we obtain
\[
\sigma_\nu^2=\sigma_\tau^2=I,
\qquad
\sigma_\nu\sigma_\tau=i\sigma_3,
\qquad
\sigma_\tau\sigma_\nu=-i\sigma_3,
\]
where $\sigma_3=
\begin{pmatrix}
1&0\\
0&-1
\end{pmatrix}
$.

The following lemma can be deduced directly from the projected Robin boundary condition satisfied by the square of the Dirac operator with infinite mass boundary conditions; see \cite[Eq.~(3.37)]{Vassilevich2003} and \cite[Remark~2.3]{EspositoGilkeyKirsten}. For convenience, we provide a direct proof.

\begin{lemma}\label{lem:boundary-log-derivative}
Let $(\lambda,\psi)$ be a solution of \eqref{eq:mit-model}. Then on $\partial\Omega$, we have
\begin{equation}\label{eq:boundary-square-norm}
\partial_\nu|\psi|^2
=
-\kappa|\psi|^2.
\end{equation}
 In particular, on $\partial\Omega\setminus Z(\psi)$, we have
\begin{equation}\label{eq:boundary-log-final}
\partial_\nu\log|\psi|
=
-\frac{\kappa}{2}.
\end{equation}
\end{lemma}

\begin{proof}
The Dirac eigenvalue equation and the boundary decomposition \eqref{eq:dirac-boundary-decomposition},
\[
-i\sigma_\nu\partial_\nu\psi
-i\sigma_\tau\partial_s\psi
=\lambda\psi,
\]
yield
\begin{equation}\label{eq:normal-derivative-spinor}
\partial_\nu\psi
=
-\sigma_\nu\sigma_\tau\partial_s\psi
+i\lambda\sigma_\nu\psi.
\end{equation}
Since $\sigma_\nu$ is Hermitian and $\lambda\in\mathbb{R}$, taking real parts gives
\begin{equation}
\frac{1}{2}\partial_\nu|\psi|^2
=
-\Re\langle\psi,
\sigma_\nu\sigma_\tau\partial_s\psi\rangle.
\label{eq:norm-normal-preliminary}
\end{equation}

By \eqref{eq:frenet-convention}, $\partial_s\sigma_\tau=-\kappa\sigma_\nu$. Differentiating the boundary condition $\sigma_\tau\psi=\psi$ tangentially, we obtain
\[
\sigma_\tau\partial_s\psi
=
\partial_s\psi+\kappa\sigma_\nu\psi.
\]
Thus, using $\sigma_\nu^2=I$, we get
\begin{align*}
\langle\psi,\sigma_\nu\sigma_\tau\partial_s\psi\rangle
&=
\langle\sigma_\nu\psi,\partial_s\psi\rangle
+\kappa|\psi|^2.
\end{align*}
On the other hand, the anticommutation relation $\sigma_\tau\sigma_\nu=-\sigma_\nu\sigma_\tau$ and the boundary condition imply
\[
\langle\psi,\sigma_\nu\sigma_\tau\partial_s\psi\rangle
=
-\langle\sigma_\nu\psi,\partial_s\psi\rangle.
\]
Comparing these two identities, we obtain
\[
\Re\langle\psi,\sigma_\nu\sigma_\tau\partial_s\psi\rangle
=
\langle\psi,\sigma_\nu\sigma_\tau\partial_s\psi\rangle
=
\frac{\kappa}{2}|\psi|^2.
\]
Substituting into \eqref{eq:norm-normal-preliminary} gives
\[
\partial_\nu|\psi|^2=-\kappa|\psi|^2.
\]
Hence, at boundary points where $\psi\neq0$,
\[
\partial_\nu\log|\psi|
=
\frac{\partial_\nu|\psi|^2}{2|\psi|^2}
=
-\frac{\kappa}{2}.
\]
\end{proof}

Let
\[
S\coloneqq |\psi|^2=|u|^2+|v|^2.
\]
Away from the common zeros of $u$ and $v$, define
\begin{equation}\label{eq:cn-defect-definition}
d_\psi\coloneqq
\frac{4|u_z\bar{v}-u\bar{v}_{z}|^2}{S^2}.
\end{equation}
The identity \eqref{eq:cn-defect-local} in the following proposition is the expression in planar complex coordinates of the two-dimensional energy--momentum identity of
Friedrich--Kim \cite{FriedrichKim2001}. To make the nonnegative defect term explicit, we provide a direct computation.

\begin{proposition}
\label{prop:cn-defect}
Let $(\lambda,\psi)$ be a solution of \eqref{eq:mit-model}. In $\Omega\setminus Z(\psi)$, we have
\begin{equation}\label{eq:cn-defect-local}
\lambda^2+2\Delta\log|\psi|=d_\psi\geqslant 0,
\end{equation}
and
\begin{equation}
4\pi\mathfrak{N}(\psi)
=
|\Omega|\lambda^2-2\pi\chi(\Omega)
-\int_\Omega d_\psi.
\label{eq:cn-defect-global}
\end{equation}
\end{proposition}

\begin{proof}
The Dirac equation is equivalent to
\[
u_{\bar{z}}=\frac{i\lambda}{2}v,
\qquad
v_z=\frac{i\lambda}{2}u.
\]
Consequently,
\[
S_z=u_z\bar{u}+v(\bar{v})_z,
\qquad
|\nabla\psi|^2
=
2\left(|u_z|^2+|v_{\bar{z}}|^2\right)
+\frac{\lambda^2}{2}S,
\]
Since $-\Delta\psi=\lambda^2\psi$, we also have
\[
\Delta S
=
2|\nabla\psi|^2-2\lambda^2S.
\]
Substituting these expressions into
\[
2\Delta\log|\psi|
=
\Delta\log S
=
\frac{\Delta S}{S}
-\frac{4|S_z|^2}{S^2},
\]
and using the identity
\[
S\left(|u_z|^2+|v_{\bar{z}}|^2\right)-|S_z|^2
=
\left|u_z\bar{v}-u(\bar{v})_z\right|^2,
\]
we obtain \eqref{eq:cn-defect-local}.

Using the previously established finiteness of the zero set in $\overline\Omega$, let
\[
\Omega_\epsilon
\coloneqq 
\Omega\setminus
\bigcup_{p\in Z(\psi)\cap\overline\Omega}
\overline{B_\epsilon(p)}.
\]
For sufficiently small $\epsilon$, the removed disks are disjoint. By \eqref{eq:interior-zero-boundary-integral}, \eqref{eq:boundary-zero-boundary-integral}, \eqref{eq:boundary-log-final} and the Gauss--Bonnet formula
\[
\int_{\partial\Omega}\kappa\dd s
=
2\pi\chi(\Omega),
\]
the divergence theorem yields
\[
\lim_{\epsilon\downarrow0}
\int_{\Omega_\epsilon}\Delta\log|\psi|\dd x
=
-\pi\chi(\Omega)-2\pi\mathfrak{N}(\psi).
\]

Integrating \eqref{eq:cn-defect-local} over $\Omega_\epsilon$ gives
\[
\int_{\Omega_\epsilon}d_\psi\dd x
=
\lambda^2|\Omega_\epsilon|
+2\int_{\Omega_\epsilon}\Delta\log|\psi|\dd x.
\]
Since $d_\psi\geqslant 0$ and $\Omega_\epsilon$ increases to $\Omega\setminus Z(\psi)$ as $\epsilon\downarrow 0$, the monotone convergence theorem, together with the finite limit above, shows that $d_\psi\in L^1(\Omega)$ and yields \eqref{eq:cn-defect-global}.
\end{proof}

\section{\texorpdfstring{The $L^1$ Lower Bound for the Defect}{The L1 Lower Bound for the Defect}}\label{sec:lower_bound_defect}
\subsection{Stereographic projection}\label{subsec:stereographic-projection}

We recall the basic differential properties of stereographic projection.

Let $S^2\subset\mathbb{R}^3$ be the unit sphere, and define the inverse stereographic projection
\[
\Psi:\mathbb{C}\longrightarrow
S^2\setminus\{(0,0,-1)\}
\]
by
\begin{equation}\label{eq:stereographic-projection}
\Psi(w)
=
\frac{1}{1+|w|^2}
\left(2\alpha,2\beta,1-\alpha^2-\beta^2\right),
\qquad w=\alpha+i\beta.
\end{equation}

A direct calculation gives
\begin{equation}\label{eq:stereographic-metric}
|\Psi_\alpha|^2
=
|\Psi_\beta|^2
=
\frac{4}{(1+|w|^2)^2},
\qquad
\Psi_\alpha\cdot\Psi_\beta=0.
\end{equation}
Thus, $\Psi$ is conformal. Moreover,
\begin{equation}\label{eq:stereographic-cross-product}
\Psi_\alpha\times\Psi_\beta
=
\frac{4}{(1+|w|^2)^2}\Psi.
\end{equation}

Let $\Omega\subset\mathbb{C}$ be a planar domain. Let $g\in C^1(\Omega;\mathbb{C})$, and write
\[
g(z)=\alpha(x,y)+i\beta(x,y),
\qquad z=x+iy.
\]
Define the sphere-valued map
\[
\mathbf{n}\coloneqq \Psi\circ g:\Omega\to S^2.
\]
By the chain rule,
\begin{equation}\label{eq:stereographic-chain-rule}
\begin{aligned}
\mathbf{n}_x
&=\alpha_x\Psi_\alpha+\beta_x\Psi_\beta,\\
\mathbf{n}_y
&=\alpha_y\Psi_\alpha+\beta_y\Psi_\beta,
\end{aligned}
\end{equation}
where $\Psi$ and its partial derivatives are evaluated at $w=g(z)$.

Combining this with \eqref{eq:stereographic-metric}, we obtain
\begin{equation}\label{eq:stereographic-inner-products}
\begin{aligned}
|\mathbf{n}_x|^2
&=
\frac{4(\alpha_x^2+\beta_x^2)}
{(1+|g|^2)^2},\\
|\mathbf{n}_y|^2
&=
\frac{4(\alpha_y^2+\beta_y^2)}
{(1+|g|^2)^2},\\
\mathbf{n}_x\cdot\mathbf{n}_y
&=
\frac{4(\alpha_x\alpha_y+\beta_x\beta_y)}
{(1+|g|^2)^2}.
\end{aligned}
\end{equation}

On the other hand, \eqref{eq:stereographic-chain-rule} and \eqref{eq:stereographic-cross-product} yield
\begin{equation}\label{eq:stereographic-n-cross}
\mathbf{n}_x\times\mathbf{n}_y
=
(\alpha_x\beta_y-\alpha_y\beta_x)
(\Psi_\alpha\times\Psi_\beta)
=
\frac{4(\alpha_x\beta_y-\alpha_y\beta_x)}
{(1+|g|^2)^2}\mathbf{n}.
\end{equation}
It follows that
\begin{equation}\label{eq:stereographic-jacobian}
\mathbf{n}\cdot(\mathbf{n}_x\times\mathbf{n}_y)
=
\frac{4(\alpha_x\beta_y-\alpha_y\beta_x)}
{(1+|g|^2)^2}.
\end{equation}

Finally, using the identities $g_z=\frac{1}{2}(g_x-ig_y)$, $g_{\bar{z}}=\frac{1}{2}(g_x+ig_y)$
and 
\[
\begin{aligned}
|g_x|^2+|g_y|^2
&=2\left(|g_z|^2+|g_{\bar{z}}|^2\right),\\
\alpha_x\beta_y-\alpha_y\beta_x
&=|g_z|^2-|g_{\bar{z}}|^2,
\end{aligned}
\]
we deduce from \eqref{eq:stereographic-inner-products} and \eqref{eq:stereographic-jacobian} that
\begin{equation}\label{eq:stereographic-energy}
|\nabla\mathbf{n}|^2
=
\frac{8\left(|g_z|^2+|g_{\bar{z}}|^2\right)}
{(1+|g|^2)^2},
\end{equation}
and
\begin{equation}\label{eq:stereographic-area-jacobian}
\mathbf{n}\cdot(\mathbf{n}_x\times\mathbf{n}_y)
=
\frac{4\left(|g_z|^2-|g_{\bar{z}}|^2\right)}
{(1+|g|^2)^2}.
\end{equation}

\subsection{Sphere Valued Maps}\label{subsec:sphere-valued-map}
Let $(u,v)^{\mathsf T}\in\mathbb{C}^2\setminus\{(0,0)\}$. Define the sphere-valued map
\begin{equation}
n\coloneqq
\frac{(2\Re(uv),\,2\Im(uv),\,|v|^2-|u|^2)}{S}
\in S^2.
\label{eq:cn-sphere-map}
\end{equation}
If $n\in W^{1,2}(B_1,S^2)$, denote
\[
e\coloneqq |\nabla n|^2=|n_x|^2+|n_y|^2,
\qquad
j\coloneqq n\cdot(n_x\times n_y).
\]
Similar sphere-valued maps associated with spinors, together with the corresponding energy and Jacobian computations, can also be found in \cite[Section~6]{AQ}. 

\begin{proposition}\label{prop:cn-sphere}
Let $\psi=(u,v)^{\mathsf T}$ be a nontrivial Dirac eigenspinor. Then the sphere-valued map $n:\Omega\setminus Z(\psi)\to S^2$ is smooth and admits a locally Lipschitz extension across the interior zeros of $\psi$. We have the identities
\begin{equation}
\label{eq:cn-sphere-identities}
e=2(\lambda^2+d_\psi),
\qquad
j=d_\psi-\lambda^2
\end{equation}
almost everywhere in $\Omega$.
\end{proposition}

\begin{proof}
By the interior zero expansion \eqref{eq:ord-zero}, near a common zero of order $m$, translated to the origin,
\[
z^{-m}(u,\bar{v})=c+O(r),
\qquad c\neq0,
\]
where the remainder has bounded first derivatives on the punctured neighborhood. Thus, $z^{-m}(u,\bar{v})$ admits a locally Lipschitz extension that remains nonzero near the origin. Since multiplication of $(u,\bar{v})$ by a nonzero complex scalar leaves $n$ unchanged, composing this extension with the smooth normalization defining $n$ yields the desired locally Lipschitz extension across the zero.

On the set where $v\neq0$, write
\[
g=\frac{u}{\bar{v}},
\qquad n=\Psi\circ g.
\]
The local computation in \cite[Section~6.2]{AQ}, applied to this coordinate, gives
\begin{equation}\label{eq:cn-g-antiholomorphic}
\frac{4|g_{\bar{z}}|^2}{(1+|g|^2)^2}
=\lambda^2,
\end{equation}
and
\begin{equation}\label{eq:cn-g-holomorphic}
\frac{4|g_z|^2}{(1+|g|^2)^2}
=d_\psi.
\end{equation}
The stereographic identities \eqref{eq:stereographic-energy} and \eqref{eq:stereographic-area-jacobian} therefore yield
\[
e=2(\lambda^2+d_\psi),
\qquad
j=d_\psi-\lambda^2.
\]
By continuity, the above identities extend to $\Omega\setminus Z(\psi)$ and hence hold almost everywhere in $\Omega$.
\end{proof}

\begin{lemma}\label{lem:cn-rigidity}
There exists $\delta> 0$, such that for every $n\in W^{1,2}(B_1,S^2)$ with
\[
e-2j=4\quad\text{a.e.},
\]
it holds
\begin{equation}\label{eq:cn-rigidity-gap}
\int_{B_1}\frac{e+2j}{4}\geqslant \delta.
\end{equation}
where $0<\delta\leqslant\pi$.
\end{lemma}

\begin{proof}
Since $|2j|\leqslant e$, denote
\[
d=\frac{e+2j}{4}\geqslant 0.
\]
We prove it by contradiction. If \eqref{eq:cn-rigidity-gap} does not hold, there exists $\{n_k\}_{k\in\mathbb{N}}$ such that
\[
e_k=2+2d_k,\qquad j_k=-1+d_k.
\]
and $\int_{B_1}d_k\to 0$. In particular, $\int_{B_1}e_k\to 2\pi$. Thus $\left\{n_k\right\}$ is uniformly bounded in $W^{1,2}(B_1,S^2)$. After passing to a subsequence, we have
\[
n_k\rightharpoonup n\quad\text{in }W^{1,2},
\qquad
n_k\to n\quad\text{in }L^2\text{ and a.e.},
\qquad |n|=1\quad\text{a.e.}
\]

We first prove the compactness of the Jacobian $j=n\cdot(n_x\times n_y)$. For component indices $a,b,c$, set
\[
M_k^{ab}
\coloneqq
\partial_xn_k^a\partial_yn_k^b
-\partial_yn_k^a\partial_xn_k^b.
\]
Since $|M_k^{ab}|\leqslant e_k$, they are uniformly bounded in $L^1$. By the identity
\[
M_k^{ab}
=
\partial_x(n_k^a\partial_yn_k^b)
-\partial_y(n_k^a\partial_xn_k^b),
\]
the weakly compactness of $\{n_k\}$ in $W^{1,2}$, the strongly compactness of $\{n_k\}$ in $L^2$, we have 
\[
M_k^{ab} \xrightarrow{\mathscr{D}^\prime(B_1)} M^{ab}\coloneqq\partial_xn^a\partial_yn^b-\partial_yn^a\partial_xn^b.
\]
$n_k - n\to 0$ in measure since $n_k - n\xrightarrow{L^2} 0$. Then $\forall \phi \in \mathscr{D}(B_1)$, $\forall \varepsilon >0$, 
\[
\begin{aligned}
\left|\int\left(n_k^c M_k^{ab} - n^c M^{ab}\right)\phi\right| 
&\leqslant C\int\left|n_k^c-n^c\right|M_k^{ab} + \int \left|M_k^{ab}-M^{ab}\right| \phi 
\\
&\leqslant C\int_{\left\{\left|n_k^c-n^c\right|>\varepsilon\right\}} e_k + \varepsilon\int_{\left\{\left|n_k^c-n^c\right|\leqslant \varepsilon\right\}} e_k + \int \left|M_k^{ab}-M^{ab}\right| \phi
\\
&\leqslant Cm\left\{\left|n_k^c-n^c\right|>\varepsilon\right\} + \varepsilon + \int \left|M_k^{ab}-M^{ab}\right| \phi.
\end{aligned}
\]
Thus
\[
j_k\xrightarrow{\mathscr{D}^\prime(B_1)} j(n).
\]
Since, in addition,  $j_k\xrightarrow{L^1}-1$, we conclude that $j(n)=-1$ almost everywhere.

Direct computation gives
\[
\left((n_k)_x-n_k\times(n_k)_y\right)^2
=e_k+2j_k=4d_k\xrightarrow{L^1} 0.
\]
Then $(n_k)_x-n_k\times(n_k)_y\xrightarrow{L^2} 0$. On another hand,
\[
(n_k)_x\xrightarrow{\mathscr{D}^\prime(B_1)} n_x, 
\quad
n_k\times(n_k)_y \xrightarrow{\mathscr{D}^\prime(B_1)} n\times n_y.
\]
It yields
\begin{equation}\label{eq:nx-n-ny}
n_x=n\times n_y \quad \text{a.e.}
\end{equation}
Then $n_x\cdot n_y=0$ and
\begin{equation}\label{eq:unit-nx-ny}
\begin{gathered}
-1 = j(n) = n\cdot (n_x\times n_y) = n\cdot \left[ \left(n\times n_y\right)\times n_y \right] = -|n_y|^2,
\\
|n_x|^2 = n_x\cdot\left(n\times n_y\right) = n\cdot \left( n_y\times n_x \right) = -j(n) = 1.
\end{gathered}
\end{equation}
Hence $n\in W^{1,\infty}_{\mathrm{loc}}\hookrightarrow C^0$. 

Restricting $n$ to a suitable neighborhood, we choose a point outside its local image as the south pole and consider the map $g\coloneqq \Psi^{-1}\circ n$. Combining \eqref{eq:stereographic-n-cross} with \eqref{eq:nx-n-ny}, we obtain $g_z=0$. Hence, $g$ is smooth and antiholomorphic. Combining \eqref{eq:stereographic-energy} with \eqref{eq:unit-nx-ny}, we obtain 
\[
\frac{4|g_{\bar z}|^2}{(1+|g|^2)^2}=1.
\]
In particular, $g_{\bar z}\neq 0$. Then
\[
0= \Delta\log
\left(\frac{|g_{\bar z}|^2}{(1+|g|^2)^2}\right)
=
-\frac{8|g_{\bar z}|^2}{(1+|g|^2)^2}
=-2,
\]
which is a contradiction. Finally, the map $n=(\cos2x,\sin2x,0)$ satisfies all the conditions and $d=1$. Hence $d\leqslant\pi$.
\end{proof}

\section{Proof of Pleijel Type Theorem}\label{sec:Pleijel-type-thm}
\begin{proposition}\label{prop:cn-packing}
There exists $\delta> 0$ such that every eigenspinor satisfies
\begin{equation}\label{eq:cn-improved}
\mathfrak{N}(\psi)
\leqslant
\left(1-\frac{\delta}{4}\right)
\frac{|\Omega|\lambda^2}{4\pi}
+O_\Omega(\lambda)-\frac{\chi(\Omega)}2
\qquad(\lambda\to\infty).
\end{equation}
The remainder term is uniform with respect to the choice of eigenspinor.
\end{proposition}

\begin{proof}
Let $B_{1/\lambda}(a)\Subset\Omega$, and denote
\[
\widetilde n(y)\coloneqq n(a+y/\lambda).
\]
By \eqref{eq:cn-sphere-identities} in Proposition~\ref{prop:cn-sphere},
\[
\widetilde e-2\widetilde j=4,
\qquad
\frac{\widetilde e+2\widetilde j}{4}
=\lambda^{-2}d_\psi(a+y/\lambda).
\]
By Lemma~\ref{lem:cn-rigidity},
\begin{equation}\label{eq:cn-one-ball}
\int_{B_{1/\lambda}(a)}d_\psi\geqslant\delta.
\end{equation}

Partition the plane into squares of side length $2/\lambda$, and retain those squares whose closures are contained in $\Omega$. The inscribed disks of these squares are pairwise disjoint. The portion of $\Omega$ not covered by the retained squares is contained in a boundary strip of width $O(1/\lambda)$, whose area is therefore
$O_\Omega(1/\lambda)$. Consequently, the number of retained squares satisfies
\[
M_\Omega(\lambda)
=\frac{|\Omega|\lambda^2}{4}
+O_\Omega(\lambda).
\]
Thus 
\[
\int_\Omega d_\psi 
\geqslant
M_\Omega(\lambda)\int_{B_{1/\lambda}(a)}d_\psi
\geqslant
\frac{\delta|\Omega|\lambda^2}{4}+ O_\Omega(\lambda).
\]
Substituting this estimate into \eqref{eq:cn-defect-global} yields the desired conclusion.
\end{proof}

The leading Weyl term for planar neutrino billiards already appears in the semiclassical analysis of Berry and Mondragon \cite[equation (74)]{BerryMondragon1987}. For completeness, we give a short variational proof of the corresponding positive spectrum asymptotics in the present setting.

\begin{lemma}\label{lem:positive-dirac-weyl}
Let $\Omega\subset\mathbb{R}^2$ be a smooth bounded domain, and let $D$ be the massless Dirac operator with infinite mass boundary conditions. Denote its positive eigenvalues, listed in nondecreasing order and counted with multiplicity, by
\[
0<\lambda_1^+\leqslant\lambda_2^+\leqslant\cdots.
\]
Then
\[
\frac{|\Omega|}{4\pi}(\lambda_k^+)^2= k+o(k),
\qquad k\to\infty.
\]
In particular, $\lambda_k^+= O_\Omega(\sqrt{k})$.
\end{lemma}

\begin{proof}
By \cite[Theorem~1.1 and Remark~2]{BFSV},
$D$ is self-adjoint and has compact resolvent.
Set 
\[
V=\operatorname{dom}D=\{\psi\in H^1(\Omega;\mathbb C^2):\sigma_\tau\psi=\psi\text{ on }\partial\Omega\}.
\]
The closed quadratic form associated with $D^2$ is
\[
q[\psi]=\|D\psi\|_{L^2(\Omega)}^2,
\qquad \operatorname{dom}q=V.
\]
Integration by parts, together with the MIT bag boundary condition, yields
\[
q[\psi]
=
\int_\Omega|\nabla\psi|^2\dd x
+\frac12\int_{\partial\Omega}\kappa|\psi|^2\dd s.
\]
Moreover,
\[
H_0^1(\Omega;\mathbb{C}^2)
\subset V
\subset H^1(\Omega;\mathbb{C}^2),
\]
and $q$ agrees with the Dirichlet energy on $H_0^1(\Omega;\mathbb{C}^2)$. By the trace inequality, $\forall \varepsilon\in(0,1)$ there exists $C_{\varepsilon,\Omega}$ such that
\[
q[\psi]\geqslant
(1-\varepsilon)\|\nabla\psi\|_{L^2(\Omega)}^2
-C_{\varepsilon,\Omega}\|\psi\|_{L^2(\Omega)}^2.
\]

Let $\gamma_j$ denote the eigenvalues of $D^2$, and let $\mu_j^{\mathrm{N}}$ and $\mu_j^{\mathrm{D}}$ denote those of the Neumann and Dirichlet Laplacians, respectively, acting on $\mathbb{C}^2$-valued functions. All eigenvalues are listed in nondecreasing order and counted with multiplicity. The Courant--Fischer min--max principle gives
\[
(1-\varepsilon)\mu_j^{\mathrm{N}}
-C_{\varepsilon,\Omega}
\leqslant \gamma_j\leqslant\mu_j^{\mathrm{D}}.
\]
The classical Weyl law for the Laplacian implies
\[
\mu_j^{\mathrm{N}}\sim\mu_j^{\mathrm{D}}
\sim\frac{2\pi}{|\Omega|}j.
\]
Letting first $j\to\infty$ and then
$\varepsilon\downarrow0$, we obtain
\[
\gamma_j\sim\frac{2\pi}{|\Omega|}j.
\]

The antilinear isometry
\[
\mathcal C(u,v)^{\mathsf T}
\coloneqq
(\bar v,\bar u)^{\mathsf T}
\]
preserves $V$ and satisfies
$D\mathcal C=-\mathcal C D$.
Consequently, the nonzero eigenvalues of $D$
occur in positive and negative pairs
with equal multiplicities.
Writing $r=\dim\ker D$, we therefore have
\[
\gamma_{r+2k-1}=\gamma_{r+2k}
=(\lambda_k^+)^2.
\]
It follows that
\[
(\lambda_k^+)^2
\sim\frac{4\pi}{|\Omega|}k,
\]
which proves the assertion.
\end{proof}

\begin{proof}[Proof of the strict asymptotic bound in Theorem \ref{thm:cn-main}]
By Lemma \ref{lem:positive-dirac-weyl},
\[
\frac{|\Omega|}{4\pi}(\lambda_k^+)^2=k+o(k),
\qquad
\lambda_k^+=O_\Omega(\sqrt{k}).
\]
Substituting these asymptotics into \eqref{eq:cn-improved} and dividing by $k$, we obtain
\[
\sup_{\psi\in E_{\lambda_k^+}\setminus\{0\}}
\frac{\mathfrak{N}(\psi)}{k}
\leqslant
\left(1-\frac{\delta}{4}\right)
\frac{|\Omega|(\lambda_k^+)^2}{4\pi k}
+O_\Omega(k^{-1/2}).
\]
Letting $k\to\infty$ yields
\[
\limsup_{k\to\infty}
\sup_{\psi\in E_{\lambda_k^+}\setminus\{0\}}
\frac{\mathfrak{N}(\psi)}{k}
\leqslant 1-\frac{\delta}{4}.
\]
Thus we may take $\eta_*=\delta/4$.

Moreover, for sufficiently large $\lambda$, the $O_\Omega(\lambda)$ remainder in \eqref{eq:cn-improved} is bounded above by
\[
\frac{\delta}{8}\frac{|\Omega|}{4\pi}\lambda^2.
\]
Consequently, there exists $\Lambda_\Omega>0$ such that, for every positive eigenvalue $\lambda\geqslant\Lambda_\Omega$ and every $\psi\in E_\lambda\setminus\{0\}$,
\[
\mathfrak{N}(\psi)
\leqslant
\left(1-\frac{\delta}{8}\right)
\frac{|\Omega|}{4\pi}\lambda^2.
\]
The threshold $\Lambda_\Omega$ depends on $\Omega$, whereas the constant $\delta$ is independent of the domain.
\end{proof}

\section{Morse Index and Spectral Index Bounds}\label{sec:morse-index}

\begin{lemma}\label{lem:cn-index-lower}
For every $\lambda>0$, we have
\[
\operatorname{ind} Q_\lambda\geqslant N_{D^2}(\lambda^2),
\]
where $N_{D^2}(\lambda^2)$ denotes the number of
eigenvalues of the Dirichlet Laplacian strictly below
$\lambda^2$, counted with multiplicity.
\end{lemma}

\begin{proof}
For $w\in H_0^1(\Omega;\mathbb C)$, we have
\[
4\|\partial_{\bar z}w\|_2^2=\|\nabla w\|_2^2,
\qquad \operatorname{tr} w=0.
\]
Consequently,
\[
Q_\lambda[w]
=\|\nabla w\|_2^2-\lambda^2\|w\|_2^2.
\]
The span of all Dirichlet eigenfunctions with eigenvalues strictly below $\lambda^2$ is a strictly negative subspace of dimension $N_{D^2}(\lambda^2)$.
\end{proof}

\begin{proposition}\label{prop:cn-index-count}
$\forall \lambda> 0$, we have
\begin{equation}\label{eq:cn-index-count}
\operatorname{ind} Q_\lambda
=\#\{j:\lambda_j^+<\lambda\},
\end{equation}
where the positive eigenvalues are counted with their complex multiplicities.
\end{proposition}

\begin{proof}
Set
\[
A[w]=4\|\partial_{\bar z}w\|_2^2,\quad
B[w]=\|\operatorname{tr} w\|_2^2,\quad
M[w]=\|w\|_2^2.
\]
Consider the family of closed quadratic forms on $E(\Omega)$ given by
\[
q_t[w]=Q_t[w]/t=A[w]/t+B[w]-tM[w].
\]
For $t>s>0$, we have
\[
q_t[w]-q_s[w]
=-\frac{t-s}{st}A[w]-(t-s)M[w]
\leqslant-(t-s)M[w].
\]
By the min--max principle, for $t>s>0$,
\[
\mu_j(t)\leqslant\mu_j(s)-(t-s).
\]
Hence each $\mu_j$ is strictly decreasing.

To prove continuity, fix $t_0>0$ and set
\[
r_0[w]
\coloneqq q_{t_0}[w]+(t_0+1)M[w]
=\frac{A[w]}{t_0}+B[w]+M[w].
\]
For $t\in[t_0/2,3t_0/2]$, we have
\[
\left|q_t[w]-q_{t_0}[w]\right|
\leqslant C_{t_0}|t-t_0|\,r_0[w],
\quad \forall w\in E(\Omega),
\]
where one may take $C_{t_0}=2/t_0+1$. Writing $\varepsilon=C_{t_0}|t-t_0|$, we obtain, for $t$ sufficiently close to $t_0$,
\[
(1-\varepsilon)r_0[w]
\leqslant q_t[w]+(t_0+1)M[w]
\leqslant (1+\varepsilon)r_0[w],
\quad \forall w\in E(\Omega).
\]
Applying the min--max principle once more yields
\[
|\mu_j(t)-\mu_j(t_0)|
\leqslant
\varepsilon
\left(\mu_j(t_0)+t_0+1\right).
\]
Thus each $\mu_j$ is continuous on $(0,\infty)$.

By a general Poincaré inequality, there exists $c>0$ such that
\[
A[w]+B[w]\geqslant cM[w].
\]
Consequently, for sufficiently small $t>0$,
\[
q_t[w]\geqslant (c-t)M[w]>0.
\]

Let $L_j$ be the span of the first $j$ Dirichlet Laplacian eigenfunctions, and let $\nu_j^{\mathrm D}$ denote the $j$th eigenvalue, counted with multiplicity. Since $L_j\subset H_0^1(\Omega)$, for every $w\in L_j$,
\[
B[w]=0,
\qquad
A[w]=\|\nabla w\|_2^2\leqslant\nu_j^{\mathrm D}M[w].
\]
Hence, by the Courant--Fischer min--max principle, $\forall t>\sqrt{\nu_j^{\mathrm D}}$,
\[
\mu_j(t)
\leqslant
\sup_{w\in L_j\setminus\{0\}}
\frac{q_t[w]}{M[w]}
\leqslant
\frac{\nu_j^{\mathrm D}}{t}-t< 0.
\]
Therefore, each $\mu_j(t)$ crosses zero exactly once.

The kernel isomorphism \eqref{eq:cn-kernel-isomorphism} identifies these zero-crossing parameters with the positive eigenvalues of the Dirac operator, preserving
multiplicities. At $t=\lambda$, the number of eigenvalue branches that have crossed zero is precisely the total multiplicity of the positive Dirac eigenvalues strictly below $\lambda$ (see the diagram below). This yields \eqref{eq:cn-index-count}.

\begin{center}
\begin{tikzpicture}[scale=0.8]
\coordinate (A) at (-2.25, 0.00);
\coordinate (B) at (7.00, 0.00);
\coordinate (C) at (0.00, -5.00);
\coordinate (D) at (0.00, 7.00);
\coordinate (E) at (0.57, 2.94);
\coordinate (F) at (4.64, -5.67);
\coordinate (G) at (0.75, 4.50);
\coordinate (H) at (4.50, -3.75);
\coordinate (I) at (1.00, 5.75);
\coordinate (J) at (5.00, -2.75);
\coordinate (K) at (1.55, 7.21);
\coordinate (L) at (5.30, -0.79);
\coordinate (M) at (4.25, 4.25);
\coordinate (N) at (4.25, -5.00);

\draw[-Stealth, name path=axis-t] (A) -- node[right, text=black, at end] {$t$} (B);
\draw[-Stealth] (C) -- node[above, text=black, at end] {$\mu$} (D);

\draw[name path=mu1] (E) -- node[above, text=black, at start] {$\mu_1(t)$} (F);
\draw[name path=mu2] (G) -- node[above, text=black, at start] {$\mu_2(t)$} (H);
\draw[name path=mu3] (I) -- node[above, text=black, at start] {$\mu_3(t)$} (J);
\draw[name path=mu4] (K) -- node[above, text=black, at start] {$\mu_4(t)$} (L);

\draw[dashed, name path=line-lambda] (M) -- node[right, text=black, at start] {$t=\lambda$} (N);

\path[name intersections={of=mu1 and axis-t, by=S}];
\path[name intersections={of=mu2 and axis-t, by=T}];
\path[name intersections={of=mu3 and axis-t, by=U}];
\path[name intersections={of=axis-t and line-lambda, by=O}];

\path[name intersections={of=mu1 and line-lambda, by=R}];
\path[name intersections={of=mu2 and line-lambda, by=Q}];
\path[name intersections={of=mu3 and line-lambda, by=P}];

\fill[black] (O) circle (4pt);
\node[above right, black] at (O) {$\lambda$};
\fill[black] (P) circle (4pt);
\fill[black] (Q) circle (4pt);
\fill[black] (R) circle (4pt);
\fill[black] (S) circle (4pt);
\fill[black] (T) circle (4pt);
\fill[black] (U) circle (4pt);
\end{tikzpicture}
\end{center}
\end{proof}

\begin{proof}[Proof of the remaining assertions in Theorem~\ref{thm:cn-main}]
By Lemma~\ref{lem:cn-index-lower} and Weyl's law for the Dirichlet Laplacian, we have
\[
\operatorname{ind} Q_\lambda
\geqslant \frac{|\Omega|}{4\pi}
\lambda^2+o(\lambda^2).
\]
On the other hand, Proposition~\ref{prop:cn-packing} yields
\[
\mathfrak{N}(\psi)
\leqslant
\left(1-\frac{\delta}{4}\right)
\frac{|\Omega|}{4\pi}
\lambda^2+O_\Omega(\lambda).
\]
Therefore there exists $L_\Omega<\infty$ such that $\forall \lambda\geqslant L_\Omega$,
\[
\mathfrak{N}(\psi)\leqslant\operatorname{ind} Q_\lambda.
\]
For $0< \lambda< L_\Omega$, the global estimate
\eqref{eq:cn-defect-global} gives
\[
\mathfrak{N}(\psi)
\leqslant
\frac{|\Omega|}{4\pi}
L_\Omega^2-\frac{\chi(\Omega)}2.
\]
Set
\[
C_\Omega\coloneqq
\max\left\{
0,\,
\frac{|\Omega|}{4\pi}
L_\Omega^2-\frac{\chi(\Omega)}2
\right\}.
\]
Combining this with $\operatorname{ind} Q_\lambda\geqslant 0$, we obtain
the estimate valid throughout the entire positive
spectrum:
\[
\mathfrak{N}(\psi)\leqslant\operatorname{ind} Q_\lambda+C_\Omega.
\]

Finally, let $\lambda=\lambda_k^+$.
We have
\[
\#\{j:\lambda_j^+<\lambda_k^+\}\leqslant k-1.
\]
The desired estimate \eqref{eq:cn-main-courant}
then follows from Proposition~\ref{prop:cn-index-count}.
\end{proof}

\section*{Declaration of AI use}
During the preparation of this manuscript, the author used OpenAI Chatgpt Plus without Codex to assist with drafting and revising mathematical arguments, checking calculations and the consistency of proofs, improving the organization and language of the manuscript, editing the \LaTeX{} source. The authors take full responsibility for the mathematical results, proofs, references, and final content of the manuscript.

\end{document}